\documentclass[12pt]{article}

\usepackage{amsmath,amssymb,latexsym, amsfonts,amsthm}
\usepackage[dvips]{graphicx}

\newtheorem{teor}{Theorem}[section]
\newtheorem{defin}[teor]{Definition}
\newtheorem{propo}[teor]{Proposition} 
\newtheorem{obs}[teor]{Remark} 
\newtheorem{theo}[teor]{Theorem}
 
\newtheorem{example}[teor]{Example} 
\newtheorem{coro}[teor]{Corollary} 

\def\proof{{\noindent \bf Proof:} \hspace{0.1 cm}}
\newcommand{\cqd}{\hfill \rule{2mm}{2mm}\vspace{.3cm}}

\def \Z {\mathbb Z}

\def \F {\mathbb F}

\begin{document}

\title{A recursive construction of units in a class of rings}

\author{ Fernanda D. de Melo-Hern\'andez \thanks{Corresponding author}
\thanks{F. D. de Melo-Hernandez (fdmelo@uem.br), Departamento de Matem\'atica, Universidade Estadual de Maring\'a, Av. Colombo 5790, 87020-900, Maring\'a, PR, Brazil.} \and
C\'esar A. Hern\'andez-Melo \thanks{C\'esar A. Hern\'andez-Melo (cahmelo@uem.br), Departamento de Matem\'atica, Universidade Estadual de Maring\'a, Av. Colombo 5790, 87020-900, Maring\'a, PR, Brazil.} 
\and Horacio Tapia-Recillas \thanks{Horacio Tapia-Recillas (htr@xanum.uam.mx), Departamento de Matem\' aticas, Universidad Aut\' onoma Metropolitana-Iztapalapa, Av. San Rafael Atlixco 186, 09340, CDMX, M\' exico.}}

\maketitle


\begin{abstract}
Let $R$ be an associative ring with identity and let $N$ be a nil ideal of $R$. It is shown that units of $R/N$ can be lifted to units in $R$. Under some mild conditions on the ring, a procedure is given to determine those lifted units in a recursive way. As an application, the units of several classes of rings are determined, including: matrix rings, chain rings, and group rings where the ring is a chain ring. Numerical examples are given illustrating the main results.



\end{abstract}


\section{Introduction}

The problem of determining the structure of the group of units of a ring has been the subject of considerable research over time. Information about the group of units of a ring can be found in the literature.  
For instance, if the ring is a finite field $\F$ its group of units is almost the entire ring: $\F^{\ast}=\F \setminus\{0\}$. If the ring is local commutative or finite commutative, an answer is provided in {\cite[Chapter XVIII pag. 353 and Chapter XXI pag. 398]{mcdonald}}, respectively. In \cite{dolz} all nonisomorphic rings whose group of units is isomorphic to a finite group of $n$ elements, where $n$ is a power of a prime or any product of prime powers not divisible by $4$, are determined.
In \cite{h-l-m} the authors determine the structure of the group of units of a finite commutative chain ring, and in \cite{k-k-s} the units of the matrix rings $M_2(\Z_b[x])$ for $b=2,3$ are determined. Other descriptions of the units of a ring can also be found in the literature, although in general it is an unresolved question. 

Of course the determination of the group of units of a ring is closely related to Fuchs' problem on determining which groups can be the groups of units of a ring. Answers to this question have been obtained in some cases. For instance, in \cite{ch-l} a class of $p$-groups is considered. In \cite{delc-d} the authors provide an answer to Fuchs' problem when the ring is torsion-free, and when the ring has characteristic zero. However, Fuchs' problem is also far from being solved in general.

The knowledge of a unit and its inverse in a ring is important in several instances, for example, in matrix theory. The group of units of a ring is also related with other areas of mathematics as well as physics, chemistry and economics. For example in \cite{k},\cite{con} a relationship with number theory and geometry, particularly with Dirichlect's Unit Theorem is considered.
  
There are several approaches to obtain the units and their inverse of a ring, one of which is to determine the units by {\it lifting} units from another ring, usually a quotient ring. For instance in \cite{perera} the question of lifting units to a ring $R$ from the ring $R/I$ for the class of separative exchange ideals $I$ and a relation with $C^{\ast}$-algebras is treated. In \cite{m-m} the following question is studied: If $R$ is a right self-injective ring and $I$ an ideal of $R$, when can a unit of $R/I$ be lifted to units of $R$? In \cite{janez} lifting units in clean rings are considered.
 
In this manuscript the units of a ring $R$ with certain properties will be determined by lifting units from a quotient ring $R/I$. This process gives a precise expression for determining the inverse of a unit of the ring $R$. Several consequences follow from this result and the units are determined in cases which include: matrix rings with entries in modular integers; commutative rings containing a nilpotent ideal; group rings $RG$ not necessarily commutative; commutative group rings $RG$ where the ring $R$ is a chain ring; and the commutative group ring $\Z_mG$ where $\Z_m$ is the ring of integers modulo $m$.

The manuscript is organized as follows: in Section 2 basic facts on units of a ring and a quotient ring are discussed. Section 3 is devoted to the recursive construction of the inverse of a unit of a ring with a certain property (see Definition 3.2), which is lifted from a unit of a quotient ring; and in Section 4, by using the main result of the previous Section, the units of the rings mentioned above are determined. Numerical examples are provided illustrating the discussed results.


\section{Basic facts}

Recall that an element $x$ of a ring $R$ is called {\it invertible} if there exists $y\in R$ such that $xy=yx=1$. The element $y$ is unique and it is denoted by $y=x^{-1}$ (the inverse of $x$). In addition, the set of units denoted by  $R^{\ast}$ of a ring $R$ form a group under the multiplication of the ring. 

The starting point for the results presented in this manuscript is the construction of inverse elements on a ring $R$ from those of a quotient ring $R/N$, where $N$ is a nil ideal. This is a particular version of the result established in {\cite[Lemma 2.1]{dolz}}, where $N$ is the Jacobson radical. 
For the sake of completeness, the proof is presented here.

\begin{propo}\label{jacobson}
Let $R$ be a ring, $N$ a nil ideal of $R$ and  \hspace{2mm}  $\bar{} : R \longrightarrow R/N$ the canonical homomorphism from $R$ to the quotient ring $R/N$. Then, $\bar f=f+N\in (R/N)^{\ast}$ if and only if $f+N\subset R^{\ast}$.  
\end{propo}

\proof
First it is shown that if $\bar f$ is an invertible element in $R/N$, then for all $x\in f+N$, $x$ is an invertible element in $R$. Since $\bar{x}=\bar{f}$ is an invertible element in the quotient ring $R/N$, there exists $\bar{g} \in R/N$ such that $xg-1 \in N$, $gx-1\in N$ and, since $N$ is a nil ideal, $(xg-1)^{2n+1}=(gx-1)^{2m+1}=0$ for some integers $n,m>0$. It is easy to see that:
\[
xe_1=1,\hspace{1.0cm}\text{with}\hspace{1.0cm}e_1=g\left[\sum_{i=0}^{2n}\binom{2n+1}{i}(xg)^{2n-i}(-1)^{i}\right],
\]
and,
\[
e_2x=1,\hspace{1.0cm}\text{with}\hspace{1.0cm}e_2=\left[\sum_{i=0}^{2m}\binom{2m+1}{i}(gx)^{2m-i}(-1)^{i}\right]g.
\]
Since $R$ is associative, it follows that $e_1=e_2$. Hence, if $x^{-1}=e_1$, $xx^{-1}=x^{-1}x=1$, i.e., $x$ is an invertible element in $R$. The proof of the other implication of the statement is straightforward.  
\cqd

Given an invertible element $\bar f\in R/N$, each element $x\in f+N\subset R$ will be called a {\it lifted invertible element} associated to $\bar f$. The previous result implies that in order to compute the inverse $x^{-1}$ of a lifted invertible element $x\in\bar{f}$, it is necessary to compute the inverse $\bar{g}$ of $\bar{x}$ in $R/N$ and then, calculate any of the expressions for $e_1=x^{-1}$ or $e_2={x^{-1}}$ given in the proof of Proposition \ref{jacobson}. One of the goals of this manuscript is to give a simplified expression to calculate $e_1=e_2=x^{-1}$. Another fundamental implication of Proposition \ref{jacobson} concerns the relation between the cardinality of the group $R^*$ and the cardinality of the group $(R/N)^*$, as noted in the following,

\begin{obs}\label{Ob1}
The following observations are easy consequences of the previous result.

\begin{enumerate}
 \item If $R$ is a ring and $N$ is a nil ideal of $R$, the group $R^{\ast}$ of invertible elements of $R$ can be described as a disjoint union of equivalence classes of invertible elements of the quotient ring $R/N$. That is,
\begin{equation}\label{fun}
R^{\ast}=\bigcup_{\bar f\in (R/N)^{\ast}}(f+N).
\end{equation}
Indeed, if $y\in R^{\ast}$ then $\bar y\in (R/N)^{\ast}$. Since $y\in\bar y$ it follows that $R^{\ast}\subset\bigcup_{\bar f\in (R/N)^{\ast}}\bar f$. The other inclusion is a consequence of Proposition \ref{jacobson}. 
\item With the same assumptions as in the previous observation, 
from relation (\ref{fun}) and Lagrange theorem: 
\[
|R^{\ast}|=\sum_{\bar f\in (R/N)^{\ast}} (f+N)=|N|\sum_{\bar f\in (R/N)^{\ast}}1=|N||(R/N)^{\ast}|.
\]
\end{enumerate} 
\end{obs}

Finally, note that Proposition \ref{jacobson} is not true if the hypothesis of $N$ being a nil ideal is removed. For instance, consider $R=\mathbb{Z}_{12}$ and $N=2\mathbb{Z}_{12}$, $1+N$ is an invertible element in $R/N$, but $3\in 1+N$ is a non-invertible element in $R$.


\section{The recursive construction}
In this section a method is provided to compute the inverse $x^{-1}$ of any lifted invertible element $x\in \bar{f}\in (R/N)^*$ given in Proposition \ref{jacobson}, under certain conditions on the ring $R$.  More precisely, if a collection $\{N_1,\dots,N_k\}$ of ideals of a ring $R$ satisfies the CNC-condition (Definition \ref{PosLiftIde}), and $x+N_1=f+N_1$ is an invertible element in the ring $R/N_1$, then $x^{-1}$ can be computed as a product of $g$ and a power of $xg$ in the ring $R$, where $g$ is any element in the class $\bar{f}^{-1}$ (Theorem \ref{IdemGeral}). For this purpose we need the following,

\begin{propo} \label{potencia}
Let $R$ be a ring and $N$ a nilpotent ideal of index $t \geq 2$ in $R$. Then the following statements holds:
\begin{enumerate}
	\item\label{pot} For any prime number $p$ such that $p \geq t$ and for all $n\in N$, there exists $r\in R$ such that
	\[
	(1+n)^p = 1 + pnr.
	\]
	\item\label{port} Let $\bar f$ be an invertible element in the quotient ring $R/N$. If a natural number $s>1$ exists such that $sN=0$, and all the prime factors of the number $s$ are greater than or equal to the nilpotency index $t$ of the ideal $N$, then $f$ is an invertible element in $R$ and  
	\[f(g(fg)^{s-1})=(fg)^{s}=1,\]
	where $g\in R$ is such that $\bar f\bar g=\bar 1$. Thus, the inverse of the invertible element $f$ is given by 
	\[f^{-1}=g(fg)^{s-1}.\]
\end{enumerate}
\end{propo}

\proof
\begin{enumerate}
\item  For $n \in R$, $n^t = 0$,
\[(1 + n)^p = \sum_{j=0}^{p}\binom{p}{j}n^j = 1 + \sum_{j=1}^{t-1}\binom{p}{j}n^j.\]
Since $p$ is a prime number, $p$ divides $\binom{p}{j}$ for all $1\leq j\leq p-1$.  Also, since $t\leq p$,
\[(1+n)^p=1+pn\left(k_1+{k_2}n+\cdots + k_{t-1}n^{t-2}\right)\]
where $k_i=\tbinom{p}{i}/p$. Therefore,
\[(1+n)^{p}= 1+pnr,\]
with $r=k_1+k_2n+\cdots + k_{t-1}n^{t-2}\in R$.

\item Since $N$ is a nilpotent ideal of $R$ and $f+N$ is an invertible element in $R/N$, Proposition \ref{jacobson} implies that $f$ is an invertible element in $R$. Let $p_1,p_2,p_3,\ldots ,p_m$ be the prime numbers, not necessarily different, appearing in the prime factor decomposition of the integer $s$, with $p_i \geq t$,
for $i=1,2,3,\cdots,m.$ Since $\bar f\bar g=\bar 1$, there exists $n \in N$ such that $fg = 1 + n$. Since $p_1\geq t$, from claim \ref{pot}, there exists $r_1\in R$ such that
\[(fg)^{p_1} = (1 + n)^{p_1}= 1 + p_1nr_1.\]
Similarly, since $p_2\geq t$ and $p_1nr_1\in N$, item \ref{pot} implies that there exists $r_2\in R$ such that 
\[(fg)^{p_1p_2}=((fg)^{p_1})^{p_2}=(1+p_1nr_1)^{p_2}=1+p_2(p_1nr_1)r_2.\]
Following this argument $r_3, r_4,\ldots, r_m\in R$ exists such that  
\[(fg)^{s}=1+sn(r_1r_2\cdots r_m).\]
In other words, 
\[(fg)^{s}=1+sh,\]
where $h=nr_1r_2\cdots r_m\in N$. Finally, since $h\in N$ and $sN=0$, $1=(fg)^{s}$ and
\[f^{-1}=g(fg)^{s-1}.\]
\end{enumerate}
\cqd

It is evident that the hypothesis in item \ref{port} of Proposition \ref{potencia}, which requires that all prime factors of $s$ be greater or equal than the nilpotency index $t$ of the ideal $N$, restricts enormously the number of applications of that proposition. For instance, if we consider $2\leq t$,  $R=\mathbb{Z}_{2^t}$ and $N=2\mathbb{Z}_{2^t}$, it is clear that $N^{t}=\{0\}$  and $2^{t-1}N=\{0\}$. Thus, according to item \ref{port} of Proposition \ref{potencia}, in order to describe the invertible elements of the ring $R$ in terms of the invertible elements of the quotient ring $R/N\cong \mathbb{Z}_2$, it is necessary that $t\leq 2$, thus $t=2$. Hence, we can only describe the invertible elements of $\mathbb{Z}_4$ in terms of the invertible elements of $\mathbb{Z}_2$, which is easily done by hand. In the following lines, we show how to overcome such restrictions.


\begin{defin}{\cite[Definition 3.2]{liftidemp}}\label{PosLiftIde}
We say that a collection $\{N_1,..., N_k\}$ of ideals of a ring $R$ satisfies the {\it CNC-condition} if the following properties hold:
\begin{enumerate}
\item \label{chc} {\bf Chain condition:} $\{0\}=N_k\subset N_{k-1}\subset\cdots \subset N_{2}\subset N_{1}\subset R$.
\item \label{nic} {\bf Nilpotency condition:} for $i=1,2,3,\ldots,k-1$, there exists $t_i \geq 2$ such that $N_i^{t_i}\subset N_{i+1}$.
\item \label{cac} {\bf Characteristic condition:} for $i=1,2,3,\ldots,k-1$, there exists $s_i\geq 1$ such that  $s_iN_i\subset N_{i+1}$. In addition, the prime factors of  $s_i$ are greater than or equal to $t_i$.   
\end{enumerate}
The minimum number $t_i$ satisfying the nilpotency condition will be called the nilpotency index of the ideal $N_i$ in the ideal $N_{i+1}$. Similarly,  
the minimum number $s_i$ satisfying the characteristic condition will be called the characteristic of the ideal $N_i$ in the ideal $N_{i+1}$. 
\end{defin}

The nilpotency condition and the characteristic condition of the previous definition can be stated as follows:

\begin{itemize}

\item[a.]\label{pronic} The nilpotency condition is equivalent to the following condition: for $i=1,2,\ldots,k-1$, $N_{i}/N_{i+1}$ is a nilpotent ideal of index $t_i$ in the ring $R/N_{i+1}$, (for details see {\cite[Definition 3.2]{liftidemp}}).  

\item[b.]\label{procac} The characteristic condition is equivalent to the following condition: for $i=1,2,\ldots,k-1$, there exists a natural number $s_i\geq1$ such that $s_i(N_{i}/N_{i+1})=0$ in the ring $R/N_{i+1}$, (for details see {\cite[Definition 3.2]{liftidemp}}).

\end{itemize}

\begin{theo}\label{IdemGeral} 
Let $	R$ be a ring, $\{N_1, N_2, \ldots, N_k\}$ a collection of ideals of $R$ satisfying the CNC-condition and $s_i$ the characteristic of the ideal $N_i$ in the ideal $N_{i+1}$.  Then, $f+N_1$ is an invertible element in the ring $(R/N_1)$ if and only if $f + N_1 \subset R^{\ast}$. Furthermore for each $x \in f+N_1$, 
\[x^{-1}=g(xg)^{s_1s_2\cdots s_{k-1}-1},\]
where $g\in R$ is such that $(f+N_1)(g+N_1)=1+N_1$. Moreover,
\begin{equation}\label{card}
|R^{\ast}|=|(R/N_1)^{\ast}||N_1|.
\end{equation}
\end{theo} 

\proof
 Since $(f+N_1)(g+N_1)=1+N_1$, for any $x \in f + N_1$, $(x+N_1)(g+N_1)=1+N_1$. The isomorphism 
\begin{equation}\label{iso1}
R/N_1\cong \frac{(R/N_2)}{(N_1/N_2)},
\end{equation}
implies that $(x+N_2+N_1/N_2)(g+N_2+N_1/N_2)=1+N_2+N_1/N_2$, so  $x+N_2+N_1/N_2$ is an invertible element in the ring $\frac{(R/N_2)}{N_1/N_2}$. Since $N_{1}/N_{2}$ is a nilpotent ideal of index $t_1$ in the ring $R/N_{2}$ and $s_1$ satisfies the hypothesis of claim 2 of Proposition \ref{potencia}, it follows that  $x+N_2$ is invertible in $R/N_2$ and
\[
(x+N_2)(g(xg)^{s_1-1}+N_2)=(xg)^{s_1}+N_2=1+N_2.
\]
From the isomorphism
\begin{equation}\label{iso2}
R/N_2\cong \frac{(R/N_3)}{(N_2/N_3)},
\end{equation}
it follows that $(x+N_3+N_2/N_3)(g(fg)^{s_1-1}+N_3+N_2/N_3)=1+N_3+N_2/N_3$, so $x+N_3+N_2/N_3$ is an invertible element in the ring $\frac{(R/N_3)}{N_2/N_3}$. Since $N_{2}/N_{3}$ is a nilpotent ideal of index $t_2$ in the ring $R/N_{3}$ and $s_2$ satisfies the hypothesis of claim 2  
of Proposition \ref{potencia}, it follows that  $x+N_3$ is invertible in $R/N_3$ and
\[
(x+N_3)(g(xg)^{s_1s_2-1}+N_3)=(xg)^{s_1s_2}+N_3=1+N_3.
\]
Since 
\begin{equation}\label{isogen}
R/N_i\cong \frac{(R/N_{i+1})}{(N_i/N_{i+1})},
\end{equation}
$x+N_{i+1}$ is invertible in $R/{N_{i+1}}$ and  
\[
(x+N_{i+1})(g(xg)^{s_1s_2\cdots s_i-1}+N_{i+1})=(xg)^{s_1s_2\cdots s_i}+N_{i+1}=1+N_{i+1}.
\]
Finally, in the last step of the chain of ideals, $x$ is invertible in $R$  and  
\[
x(g(xg)^{s_1s_2\cdots s_{k-1}-1})=(xg)^{s_1s_2\cdots s_{k-1}}=1,
\]
which implies that  $f + N_1 \subset R^{\ast}$ and $x^{-1}=g(xg)^{s_1s_2\cdots s_{k-1}-1}$. Claim (\ref{card}) follows from the isomorphism (\ref{isogen}) and relation (\ref{fun}),
\begin{equation}\label{ok1}
|(R/N_{i+1})^{\ast}|=|(R/N_i)^{\ast}||N_i/N_{i+1}|,
\end{equation}
for $i=1,2,\dots, k-1$. Applying the previous relation ($k-1$) times, 
\[
|R^{\ast}|=|(R/N_{k})^{\ast}|=|(R/N_{1})^{\ast}||N_{1}/N_{2}||N_{2}/N_{3}|\cdots |N_{k-2}/N_{k-1}||N_{k-1}/N_{k}|.
\]
Finally, from Lagrange's theorem,
\[
|R^{\ast}|=|(R/N_{1})^{\ast}||N_{1}|/|N_{k}|=|(R/N_{1})^{\ast}||N_{1}|,
\]
proving the result.
\cqd

\begin{obs}
Observe that if $\{N_1,N_2,\ldots,N_k\}$ is a collection of ideals of the ring $R$ satisfying the CNC-condition, any invertible element $x+N_1$ of the ring $R/N_1$ is lifted up to the invertible element $x+N_2$ of the ring $R/N_2$. This new invertible element is lifted up to the invertible element $x+N_3$ of the ring $R/N_3$, and so on. At the end of this process, $x$ is an invertible element in $R$ with $x^{-1}=g(xg)^{s_1s_2\cdots s_{k-1}-1}$. The following chain of ring homomorphisms, 
\[
R\xrightarrow{\phi_{k-1}}\frac{R}{N_{k-1}}\xrightarrow{\phi_{k-2}} \cdots \xrightarrow{\phi_{3}}\frac{R}{N_{3}}\xrightarrow{\phi_{2}}\frac{R}{N_{2}}\xrightarrow{\phi_{1}}\frac{R}{N_{1}},
\]
appears naturally in the lifting process of the invertible element $x+N_1 \in R/N_1$. 
\end{obs}

\begin{obs}
If in Theorem \ref{IdemGeral} it is assumed that the ring $R$ is commutative, the inverse $x^{-1}$ of $x$ in the ring $R$ can be given as
\[
x^{-1}=g^{s_1s_2\cdots s_{k-1}}x^{s_1s_2\cdots s_{k-1}-1},
\]
where $(x+N_1)(g+N_1)=(1+N_1)$. If $(x+N_1)(g+N_1)=(1+N_1)$ and $g+N_1=h+N_1$, then 
\[
x^{-1}=h^{s_1s_2\cdots s_{k-1}}x^{s_1s_2\cdots s_{k-1}-1}.
\]
Since $x^{-1}$ is unique, $g^{s_1s_2\cdots s_{k-1}}=h^{s_1s_2\cdots s_{k-1}}$. Consequently, the function 
\[
\psi:(R/N_1)^{\ast}\rightarrow R^{\ast},\hspace{0.3cm} g+N_1\rightarrow \psi(g+N_1)=g^{s_1s_2\cdots s_{k-1}},
\]
defines a group homomorphism. Furthermore, the classical isomorphism $I:R^{\ast}\rightarrow R^{\ast}$ that maps $x\in R^{\ast}$ to its inverse $x^{-1}$, takes the following explicit expression:
\[
I(x)=\psi\left((x+N_1)^{-1}\right)x^{s_1s_2\cdots s_{k-1}-1}.
\]
Hence, the inverse of an invertible element $x$ can be computed as a multiplication of a power of $x$ by a power of an element in the equivalence class of the inverse of $x+N_1\in R/N_1$. 
\end{obs}


\section{Applications of the main result}

\quad In this section Theorem \ref{IdemGeral} will be applied in order to determine the group of units of several rings which include: rings containing a nilpotent ideal; matrix ring $M_n(R)$ where $R$ is a commutative ring containing a collection of ideals satisfying the CNC-condition; the group ring $RG$ where the ring $R$ contains a nilpotent ideal; group rings $RG$ where $R$ is a chain ring; and the group ring $\Z_mG$ where $\Z_m$ is the ring of integers modulo $m$. Examples are given illustrating the results.



\subsection{Rings containing a nilpotent ideal}

In the following result, for a ring $R$ containing a nilpotent ideal $N$, starting from the set of invertible elements of the quotient ring $R/N$, the set of invertible elements of $R$ is determined. 

\begin{propo}\label{GeralNil} 
Let $R$ be a ring and $N$ a nilpotent ideal of nilpotency index $k\geq 2$ in $R$.  Let $s>1$ be the characteristic of the quotient ring $R/N$. Then $f+ N$ is an invertible element in $R/N$ if and only if $f + N \subset R^{\ast}$. Furthermore for each $x \in f + N$,
\[x^{-1}=g(xg)^{(s^{k-1}-1)}\]
where $g\in R$ is such that $\bar{f}\bar{g}=\bar{1}$. Moreover, $|R^{\ast}|=|(R/N)^{\ast}||N|$.
\end{propo} 
\proof The proof of this proposition is a consequence of Theorem \ref{IdemGeral}. First it is shown that the collection  $B=\{N, N^2,...,N^k\}$ of ideals of the ring $R$ satisfies the CNC-condition with nilpotency index and characteristic of the ideal $N^{i}$ in the ideal $N^{i+1}$ being $t_i=2$ and $s_i=s$ for all $i=1,2,3,\dots,k-1$. Thus, 
\begin{enumerate}
\item It is clear that the collection $B$ satisfies the chain condition. 
\item Since $(N^i)^2 = N^{2i}$ and $i+1\leq 2i$ for all $i=1,2,3,\dots, k-1$, it follows that $(N^i)^2 \subset N^{i+1}$. Hence, the collection $B$ satisfies the nilpotency condition. 
\item Since the ring $R/N$ has characteristic $s$, there exists $n\in N$ such that $\sum_{i=1}^s1_R=n$. Since 
\begin{equation}
sN^{i} =(1_R+\cdots+1_R)N^{i}=nN^{i}\subset N^{i+1}, 
\end{equation}
it follows that $sN^i\subset N^{i+1}$ for all $i=1,2,3,\dots,k-1$. In addition, all prime factors of $s_i=s$ are greater or equal to the nilpotency index $t_i=2$, proving that the collection $B$ satisfies the characteristic condition.
\end{enumerate}
Theorem \ref{IdemGeral} implies that $f + N \subset R^{\ast}$ and the inverse of an invertible element  $x \in f + N$  is
$x^{-1}=g(xg)^{(s^{k-1}-1)}$. 
\cqd

An application of the previous result is the following:

\begin{coro}\label{coroe}
Let $R$ be a commutative ring, $a \in R$ be a nilpotent element of index $k$ and $s>1$ be the characteristic of the quotient ring $R/\langle a\rangle$. Then, $f+ \langle a \rangle$ is an invertible element in $R/\langle a\rangle$ if and only if $f + \langle a \rangle \subset R^{\ast}$. Furthermore for each $x \in f + \langle a \rangle$,
\[x^{-1}=g(xg)^{(s^{k-1}-1)}\]
where $g\in R$ is such that $\bar{f}\bar{g}=\bar{1}$. Moreover, $|R^{\ast}|=|(R/\langle a\rangle)^{\ast}||\langle a\rangle|.$
\end{coro} 
\proof Since $R$ is a commutative ring, $\langle a \rangle$ is a nilpotent ideal of nilpotency index $k$ in $R$, and the result follows immediately from Proposition \ref{GeralNil}
\medskip
\cqd

\begin{obs}\label{nilcom}
If in Proposition \ref{GeralNil} it is assumed that the ring $R$ is commutative, the inverse $f^{-1}$ of $f\in R$ can be written as
\[
f^{-1}=g^{s^{k-1}}f^{s^{k-1}-1},
\]
where $g\in R$ is such that $\bar{f}\bar{g}=\bar{1}$.
\end{obs}

\noindent

\begin{example}
Let 
\[
\mathbb{Z}_{p^{k}}[i]=\{a+bi: a,b\in \mathbb{Z}_{p^k}, i^2=-1\},
\] 
where $p>2$ is a prime and $k>1$ a natural number. 
It is easy to see that $a=p$ is a nilpotent element of index $k$ in the ring $\mathbb{Z}_{p^{k}}[i]$. Since 
\[
\frac{\mathbb{Z}_{p^{k}}[i]}{\langle p\rangle}\cong \mathbb{Z}_{p}[i],
\]
and the ring $\mathbb{Z}_{p}[i]$ has characteristic $s=p$, if $ f + \langle p \rangle$ is a unit in $\frac{\mathbb{Z}_{p^{k}}[i]}{\langle p\rangle}$, Corollary \ref{coroe} and Remark \ref{nilcom} imply that if $x\in f + \langle p \rangle $ then 
\[
x^{-1}=g^{s^{k-1}}x^{s^{k-1}-1},
\]
and $|\mathbb{Z}_{p^{k}}[i]^{\ast}|=|(\mathbb{Z}_{p^{k}}[i]/\langle p\rangle)^{\ast}||\langle p\rangle|.$

It is easy to see that for any even prime $\mathbb{Z}_{p}[i]$ is a finite field with $p^2$ elements. Hence  
$(\mathbb{Z}_{p}[i])^{\ast}$ is a cyclic group with $(p^2-1)$ elements. Therefore, 
\[
(\mathbb{Z}_{p^{k}}[i])^{\ast}=\bigcup_{u \in (\mathbb{Z}_{p}[i])^{\ast}} (u+ \langle p \rangle),
\]
with a total of $(p^2-1)p^{k}$ elements.

For example, if $p=3$ and $k=2$, it is not difficult to see that $(\mathbb{Z}_{3}[i])^{\ast} = \{1,2, i, 2i, 1 + i, 2 + i, 1 +2i, 2+2i\}.$ Hence,
\[
(\mathbb{Z}_{3^{2}}[i])^{\ast}=\bigcup_{u \in (\mathbb{Z}_{3}[i])^{\ast}} (u+ \langle 3 \rangle),
\]
with a total of $(3^2 -1)\times 3^2= 72$ elements.
\end{example}


\subsection{Matrix ring}

The following result describes the group of units of the matrix ring $M_n(R)$ obtained from the units of the matrix ring $M_n(R/N_1)$, where $R$ is a commutative ring containing a collection of ideals $\{N_1,\dots,N_k\}$ satisfying the CNC-condition.  
\begin{propo}\label{MatrixRing} 
Let $R$ be a commutative ring and let $M_n(R)$ be the ring of  $n\times n$ matrices with entries from $R$. Let $\{N_1, N_2, \ldots, N_k\}$ be a collection of ideals of $R$ satisfying the CNC-condition, and  
let  $s_i$ be the characteristic of the ideal $N_i$ in the ideal $N_{i+1}$. Let $f=[f_{ij}]\in M_n(R)$, then $\overline{f}=[f_{ij}+N_1]$ is an invertible element in the matrix ring $M_n(R/N_1)$ if and only if $f + M_n(N_1) \subset (M_n(R))^\ast$. Furthermore for each $x\in f + M_n(N_1)$,
\[
x^{-1}=g(xg)^{s_1s_2\cdots s_{k-1}-1},
\] 
where $g=[g_{ij}]\in M_n(R)$ is such that $\overline{f}\overline{g}=\overline{I}$, with $I$ being the identity matrix. Moreover, 
\begin{equation}\label{cardmatrix}
|(M_n(R))^{\ast}|=|(M_n(R/N_1))^{\ast}||N_1|^{n^2}.
\end{equation}
\end{propo} 
\proof The proof of this proposition is a consequence of Theorem \ref{IdemGeral}. First, we claim that the collection  $B=\{M_n(N_1), M_n(N_2),...,M_n(N_k)\}$ of ideals of the ring $M_n(R)$ satisfies the CNC-condition with nilpotency index and characteristic of the ideal $M_n(N_i)$ in the ideal $M_n(N_{i+1})$ being exactly the same nilpotency index and characteristic of the ideal $N_i$ in the ideal $N_{i+1}$. Indeed, 
\begin{enumerate}
\item It is obvious that the collection $B$ satisfies the chain condition. 
\item If $t_i$ denotes the nilpotency index of the ideal $N_i$ in the ideal $N_{i+1}$, then $N_{i}^{t_i}\subset N_{i+1}$. It is easy to see that  $(M_n(N_i))^{t_i}\subset M_n(N_i^{t_i}))$. Hence,
\[
(M_n(N_i))^{t_i}\subset (M_n(N_{i+1})),
\]
proving that the collection $B$ satisfies the nilpotency condition.
\item If $s_i$ is the characteristic of the ideal $N_i$ in the ideal $N_{i+1}$, then $s_iN_{i}\subset N_{i+1}$. It is clear that  $s_i(M_n(N_i))=M_n(s_iN_i)$. Hence,
\[
s_i(M_n(N_i))\subset M_n(N_{i+1}).
\]
Now, since the collection $\{N_1, N_2, \ldots, N_k\}$ satisfies the CNC-condition, all prime factors of the characteristic $s_i$ are greater than or equal to the nilpotency index $t_i$ for all $i=1,2,3,\dots,k-1$. Thus the collection $B$ satisfies the characteristic condition.
\end{enumerate}
From the isomorphism 
\begin{equation}\label{isomatr}
\frac{M_n(R)}{M_n(N_1)}\cong M_n(R/N_1),
\end{equation}
it follows that $\overline{f}=[f_{ij}+N_1]$ is an invertible element in the matrix ring $M_n(R/N_1)$ if and only if $f+M_n(N_1)$ is an invertible element in the ring  $M_n(R)/M_n(N_1)$. From Theorem \ref{IdemGeral}, it follows that $f+M_n(N_1)\subset (M_n(R))^{\ast}$ and for each $x\in f + M_n(N_1)$,
\[
x^{-1}=g(xg)^{s_1s_2\cdots s_{k-1}-1},
\] 
where $g=[g_{ij}]\in M_n(R)$ is such that $\overline{f}\overline{g}=\overline{I}$. Finally, from relation (\ref{card}) in Theorem \ref{IdemGeral} and the isomorphism given in (\ref{isomatr}), we conclude that
\[
|(M_n(R))^{\ast}|=|(M_n(R/N_1))^{\ast}| |M_n(N_1)|,
\]
and the proposition is proved.
\cqd

The following results describe the set of invertible elements of the matrix ring $M_n(R)$ in terms of the matrix ring $M_n(R/N)$, where $R$ is a commutative ring containing a nilpotent ideal $N$. 

\begin{coro}\label{GeralNilMnR} 
Let $R$ be a commutative ring and let $M_n(R)$ be the ring of  $n\times n$ matrices with entries from $R$. Let $N$ be a nilpotent ideal of $R$ of index $k$ in $R$, and  let  $s$ be the characteristic of the quotient ring $R/N$. Let $f=[f_{ij}]\in M_n(R)$, then $\overline{f}=[f_{ij}+N]$ is an invertible element of the matrix ring $M_n(R/N)$ if and only if $f + M_n(N) \subset (M_n(R))^\ast$. Furthermore for each $x\in f + M_n(N)$,
\[
x^{-1}=g(xg)^{s^{k-1}-1},
\] 
where $g=[g_{ij}]\in M_n(R)$ is such that $\overline{f}\overline{g}=\overline{I}$, with $I$ being the identity matrix. Moreover, 
\begin{equation}\label{cardmatrix}
|(M_n(R))^{\ast}|=|(M_n(R/N))^{\ast}||N|^{n^2}.
\end{equation}
\end{coro}

\proof 
It is obvious that the collection $\{N, N^2,...,N^{k}\}$ of ideals of the ring $R$ satisfies the CNC-condition with constant characteristic $s_i=s$, for all $i=1,2,3,\cdots, k-1$. Thus, the result follows from Proposition \ref{MatrixRing}.
\cqd 

\begin{coro}\label{GeralEleNilMnR} 
Let $R$ be a commutative ring and $M_n(R)$ be the ring of  $n\times n$ matrices with entries from $R$. Let $a$ be a nilpotent element  of  index $k$ in $R$ and  $s$ be the characteristic of the quotient ring $R/\langle a\rangle$. Let $f=[f_{ij}]\in M_n(R)$, then $\overline{f}=[f_{ij}+\langle a\rangle]$ is an invertible element of the matrix ring $M_n(R/\langle a\rangle)$ if and only if $f + M_n(\langle a\rangle) \subset (M_n(R))^\ast$. Furthermore for each $x\in f + M_n(\langle a\rangle)$
\[
x^{-1}=g(xg)^{s^{k-1}-1},
\] 
where $g=[g_{ij}]\in M_n(R)$ is such that $\overline{f}\overline{g}=\overline{I}$, with $I$ being the identity matrix. Moreover, 
\begin{equation}\label{cardmatrix}
|(M_n(R))^{\ast}|=|(M_n(R/\langle a\rangle))^{\ast}||\langle a\rangle|^{n^2}.
\end{equation}
\end{coro}

\proof
 It is enough to observe that the ideal $N=\langle a\rangle$ is nilpotent of index $k$ in the ring $R$. The result follows from Corollary \ref{GeralNilMnR}.
\cqd

In the following lines the above results will be illustrated in the context of the ring of matrices with entries in the ring $R=\Z_{p^k}$. First, recall that a matrix $f\in M_n(\Z_{p^k})$ is invertible if and only if $\text{det}(f)$ is an invertible element of the ring $\Z_{p^k}$, where  $\text{det}(f)$ denotes the determinant of the matrix $f$. In this case the inverse of $f$ is given by the relation

\begin{equation}\label{classic}
f^{-1}=\text{det}(f)^{-1}\,\text{adj}(f),
\end{equation}
where $\text{adj}(f)$ is the adjoint matrix of $f$, that is, the transpose of its cofactor matrix. Since $p$ is a nilpotent element of index $k$ in the ring $\Z_{p^k}$ and $\Z_{p^k}/\langle p \rangle\cong \Z_{p}$ has characteristic $p$, Corollary \ref{GeralEleNilMnR} can be used to determine the inverse of the matrix $f$, if it exists.  The following steps define an algorithm to compute the inverse of a matrix $f=[f_{ij}]\in M_n(\Z_{p^k})$ by using this Corollary:  

\begin{enumerate}\label{StepsM}
\item\label{ci} Check if $f$ is invertible by recalling that $f=[f_{ij}]$ is invertible in $M_n(\Z_{p^k})$ if and only if $\overline{f}=[f_{ij}\, \text{mod}\,(p)]$ is invertible in $M_n(\Z_{p})$, which is equivalent to $\text{det}(\overline{f})\neq 0$ in $\Z_p$. 

\item\label{cip} If $\overline{f}$ is invertible in $M_n(\Z_{p})$ proceed to compute its inverse $\overline{g}$, that is, $\overline{g}\overline{f}=\overline{f}\overline{g}=\overline{I}$, by using relation (\ref{classic}) or any other known procedure for this purpose.

\item\label{cia} Finally, compute the inverse of $f\in M_n(\Z_{p^k})$ by using the relation: 
\begin{equation}\label{nonclassic}
f^{-1}=g(fg)^{p^{k-1}-1}.
\end{equation}
\end{enumerate}  

There are some advantages to using relation (\ref{nonclassic}) for computing the inverse of a matrix $f\in M_n(\Z_{p^k})$. If $x\in M_n(\Z_{p^{k}})$ is such that $\overline{x}=\overline{f}$ and $\overline{f}$ is invertible, then $x$ is invertible and the inverse of $x$ can be computed as
$
x^{-1}=g(xg)^{p^{k-1}-1},
$ where $\overline{g}$ denotes the inverse of $\overline{f}$. Hence, when $\overline{x}=\overline{f}$, $f$ is invertible and the inverse $\overline{g}$ of $\overline{f}$ has been computed in order to compute $x^{-1}$, so steps \ref{ci} and \ref{cip} in the algorithm described above can be avoided.\\

As mentioned above, the invertibility of a matrix $f\in M_n(\Z_{p^k})$ can be analyzed directly by computing $\text{det}(f)$ in the ring $\Z_{p^k}$,  or indirectly by computing  $\text{det}(\overline{f})$ in the field $\Z_{p}$. Since the number of entries equal to zero in the matrix $\overline{f}$ could be greater than the number of entries equal to zero of the matrix $f$, the computation of $\text{det}(\overline{f})$ could be easier than the computation of $\text{det}(f)$.
\begin{example}
Let $f\in M_3(\Z_{3^3})$ be given by
\begin{equation}\label{deff}
f=
\begin{bmatrix}
19 & 12 & 22 \\
6 & 5 & 24 \\
0 & 16 & 11 
\end{bmatrix}\hspace{1.0cm}\text{then,}\hspace{1.0cm}
\overline{f}=
\begin{bmatrix}
1 & 0 & 1 \\
0 & 2 & 0 \\
0 & 1 & 2 
\end{bmatrix}.
\end{equation}
It is easy to see that $\text{det}(\overline{f})=1$,  so $\overline {f}$ is invertible in $M_3(\Z_{3})$ and likewise, it is easy to check that  $\overline{g}$, the inverse of $\overline {f}$  is given by: 
\begin{equation}\label{defg}
\overline{g}=
\begin{bmatrix}
1 & 1 & 1 \\
0 & 2 & 0 \\
0 & 2 & 2 
\end{bmatrix}.
\end{equation}
In order to determine $f^{-1}$, using relation (\ref{nonclassic}) with $p=3$, $k=3$, $f$ given in (\ref{deff}) and $g=\overline{g}$ given in (\ref{defg}), we find that 
\[ \label{calinv}
f^{-1}=g(fg)^8=g
\begin{bmatrix}
19 & 6 & 9 \\
6 & 10 & 0 \\
0 & 0 & 22 
\end{bmatrix}^8=g
\begin{bmatrix}
19 & 21 & 18 \\
21 & 1 & 0 \\
0 & 0 & 16 
\end{bmatrix}
=\begin{bmatrix}
13 & 22 & 7 \\
15 & 2 & 0 \\
15 & 2 & 5 
\end{bmatrix},
\]
finishing the example.
\end{example}
Next, a way to determine the inverse of an invertible matrix in $M_n(\Z_m)$ is given. 
Let $m=p_1^{k_1}p_2^{k_2}\cdots p_j^{k_j}$ be a positive integer with its prime factorization. It is obvious that the function 
\begin{equation}\label{isomm}
\begin{aligned}
&\Psi: M_n(\Z_m)\longrightarrow M_n(\Z_{p_1^{k_1}})\times\cdots\times M_n(\Z_{p_j^{k_j}}),\\ 
&\Psi(f)=(f_1,\dots, f_j),\hspace{0.2cm}f_i=f\,\text{mod}\left(p_i^{k_i}\right) 
\end{aligned}
\end{equation}
induced by the function given by the Chinese Remainder Theorem (CRT) on $\Z_m$, is a ring isomorphism. For $m_i=m/p_i^{k_i}$ and $s_i$ such that $s_im_i=1\,\text{mod}\left(p_i^{r_i}\right)$, the inverse of $\Psi$ is given by
\[
\Psi^{-1}(h_1,\dots,h_j)=s_im_1h_1+\cdots+s_jm_jh_j.
\]
Moreover, a matrix $f$ is invertible in $M_n(\Z_m)$ if and only if each of the coordinates $f_i$ of $f$ is invertible in $M_n(\Z_{p_i^{k^i}})$. Finally, if $f$ is invertible, the inverse of $f$ can be computed by using the relation
\[\label{forfin}
f^{-1}=s_1m_1f_1^{-1}+\cdots+s_jm_jf_j^{-1},
\]
where $f_i^{-1}=g_i(f_ig_i)^{(p_i^{k_i-1}-1)}$ and $\overline {g_i}$ is the inverse of the matrix $\overline{f_i}$ in the ring $M_n(\Z_{p_i})$.


\subsection{Group rings}

If $R$ is a ring containing a collection of ideals satisfying the CNC-condition and $G$ is a group, by invoking Theorem \ref{IdemGeral}, the inverse of the units of the group ring $RG$ are determined as follows:. 

\begin{propo}\label{GeralGR} 
Let $	R$ be a ring and $G$ a group. Let $\{N_1, N_2, \ldots, N_k\}$ be a collection of ideals of $R$ satisfying the CNC-condition with $s_i$ being the characteristic of the ideal $N_i$ in the ideal $N_{i+1}$. Then, $f+N_1G$ is an invertible element in the group ring $(R/N_1)G$ if and only if $f + N_1G \subset (RG)^\ast$. Furthermore for $x \in f + N_1G$,
\[
x^{-1}=g(xg)^{s_1s_2\cdots s_{k-1}-1},
\] 
where $g\in RG$ is such that $(f+N_1G)(g+N_1G)=1+N_1G$. Moreover, $|(RG)^{\ast}|=|N_1|^{|G|}|((R/N_1)G)^{\ast}|$.
\end{propo} 

\proof 
We claim that the collection $B=\{N_1G, N_2G, \ldots, N_kG\}$ of ideals of the group ring $RG$ satisfies the CNC-condition with nilpotency index and characteristic of the ideal $N_iG$ in the ideal $N_{i+1}G$ being the same nilpotency index and characteristic of the ideal $N_i$ in the ideal $N_{i+1}$. Indeed, 
\begin{enumerate}
\item It is clear that the collection $B$ satisfies the chain condition. 
\item If $t_i$ denotes the nilpotency index of the ideal $N_i$ in the ideal $N_{i+1}$, then $N_{i}^{t_i}\subset N_{i+1}$. It is not difficult to see that  $(N_iG)^{t_i}=N_i^{t_i}G$. Hence,
\[
(N_iG)^{t_i}=N_i^{t_i}G\subset N_{i+1}G,
\]
proving that the collection $B$ satisfies the nilpotency condition.
\item If $s_i$ is the characteristic of the ideal $N_i$ in the ideal $N_{i+1}$, then $s_iN_{i}\subset N_{i+1}$. It is clear that  $s_i(N_iG)=(s_iN_i)G$. Hence,
\[
s_i(N_iG)=(s_iN_i)G\subset N_{i+1}G.
\]
Since the collection $\{N_1,N_2,\dots,N_k\}$ satisfies the CNC-condition, it is obvious that all prime factors of the characteristic $s_i$ are greater than or equal to the nilpotency index $t_i$ for all $i=1,2,3,\dots,k-1$. Thus the collection $B$ satisfies the characteristic condition.
\end{enumerate}
From Theorem \ref{IdemGeral} and the following isomorphism,
\[
\frac{RG}{N_1G}\cong \left(\frac{R}{N_1}\right)G,
\]
it follows that  $f + N_1G \subset (RG)^{\ast}$, the inverse of an  invertible element  $x \in f + N_1G$  is $x^{-1}=g(xg)^{s_1s_2\cdots s_{k-1}-1}$, and $|(RG)^{\ast}|=|N_1|^{|G|}|((R/N_1)G)^{\ast}|$.
\cqd

\begin{obs}
If it is assumed in Proposition \ref{GeralGR} that the ring $R$ and the group $G$ are commutative, the inverse $f^{-1}$ of $f\in RG$ can be expressed as
\[
f^{-1}=g^{s_1s_2\cdots s_{k-1}}f^{(s_1s_2\cdots s_{k-1})-1},
\]
where $g\in RG$ is such that $(f+N_1G)(g+N_1G)=1+N_1G$.
\end{obs}

\begin{coro}\label{GeralNilGR} 
Let $	R$ be a ring, $N$ a nilpotent ideal of index $k$ in $R$, $G$ a group and $s$ the characteristic of the quotient ring $R/N$. Then, $f+NG$ is an invertible element in the group ring $(R/N)G$, if and only if $f + NG \subset (RG)^{\ast}$. Furthermore for each $x \in f + NG,$
\[
x^{-1}=g(xg)^{s^{k-1}-1},
\] 
where $g\in RG$ is such that $(f+NG)(g+NG)=1+NG$. Moreover, $|(RG)^{\ast}|=|N|^{|G|}|((R/N)G)^{\ast}|$.
\end{coro}
\proof 
It readily seen that the collection $\{N, N^2,...,N^{k}\}$ of ideals of the ring $R$ satisfies the CNC-condition with constant characteristic $s_i=s$ for all $i=1,2,3,\cdots, k-1$. The result follows from Proposition \ref{GeralGR}.
\cqd 

\begin{example}
Consider the group ring $M_n(\Z_9) S_3$, where $S_3$ is the symmetric group in three symbols. 
Observe that $N = \langle 3 \rangle$ is a nilpotent ideal of index $2$ in the ring $\Z_9$, so the collection $\{M_n(N)S_3, M_n(N^2)S_3\}$ of ideals of the ring $M_n(\Z_9)S_3$ satisfies the CNC-codition with constant characteristic $s = 3$. Now, for $\sigma=(123) \in S_3$ and $I$ the identity matrix in $M_n(\Z_9)$, the element $ f = 2I +8 I\sigma$ is an invertible element in the group ring $M_n(\Z_9)S_3$ and $\bar{g} = I +2I\sigma + I\sigma^2 \in M_n(\Z_9/N)S_3$ is such that $\bar{f}\bar{g} = 1$ in $M_n(\Z_9/N)S_3$. By Corollary \ref{GeralNilGR}, it follows that $f + M_n(N)S_3 \subset (M_n(\Z_9)S_3)^{\ast}$ and $x^{-1} = g(xg)^2$, where $x$ denotes any element in $f + M_n(N)S_3$. For instance, since $f \in f + M_n(N)S_3$, the inverse of $f$ in the group ring $M_n(\Z_9)S_3$ is given by
\[f^{-1} = g(fg)^2 =  7I +8I \sigma + 4I\sigma^2.\]
In the same fashion,
since $x = 2I + 2I \sigma \in f + M_n(N)S_3$, then 
\[ x^{-1} = g(xg)^2 = 4I + I\sigma + 4 I\sigma^2,\]
is the inverse of $x$ in the group ring $M_n(\Z_9)S_3$.
\end{example}

\begin{coro}\label{cyu}
Let $	R$ be a commutative ring, $a$ be a nilpotent element of index $k$ in $R$, $G$ a group and $s$ the characteristic of the quotient ring $R/\langle a\rangle$. Then, $f+\langle a\rangle G$ is an invertible element in the group ring $(R/\langle a\rangle)G$ if and only if $f+\langle a\rangle G \subset ((R/\langle a\rangle)G)^{\ast}$. Furthermore for each $x \in f+\langle a\rangle G$,
\[
x^{-1}=g(xg)^{s^{k-1}-1},
\]
where $g\in RG$ is such that $(f+\langle a\rangle G)(g+\langle a\rangle G)=1+\langle a\rangle G$. Moreover, $|(RG)^{\ast}|=|\langle a\rangle|^{|G|}|((R/\langle a\rangle)G)^{\ast}|$.
\end{coro} 
\proof 
Since $R$ is a commutative ring, the ideal $N = \langle a \rangle $ is a nilpotent ideal of index $k$ in $R$, and the result follows from Corollary \ref{GeralNilGR}.
\cqd

\begin{obs}
If it is assumed in Corollaries \ref{GeralNilGR} and \ref{cyu} that the ring $R$ and the group $G$ are both commutative, the inverse $f^{-1}$ of $f\in RG$ can be expressed as:
\[
f^{-1}=g^{s^{k-1}}f^{s^{k-1}-1},
\]
where $g\in RG$ is such that $\bar f\bar g=\bar 1$.
\end{obs}


\subsection{Commutative group rings $RG$ with $R$ a chain ring}

Let $R$ be a finite commutative chain ring and $G$ a commutative group. It is well known that $R$ contains a unique maximal nilpotent ideal $N=\langle a\rangle$ for some $a\in R$. If $k$ denotes the nilpotency index of $a$, and $p$ denotes the characteristic of the residue field $\F=R/\langle a\rangle$, from Corollary \ref{cyu}, it follows that  
\begin{equation}\label{few}
(RG)^{\ast}=\left\{g^{r}f^{r-1}: r=p^{k-1}, \bar{f},\bar{g}\in (\F G)^{\ast} \hspace{0.1cm}\text{and}\hspace{0.1cm} \bar{f}\bar{g}=\bar{1}\right\},
\end{equation}
and
\begin{equation}\label{fe}
|(RG)^{\ast}|=|(\F G)^{\ast}||\langle a\rangle|^{|G|}.
\end{equation}
Examples of finite commutative chain rings include the ring of modular integers $R=\mathbb{Z}_{p^k}$, where $p$ is a prime number and $k>1$ is an integer. In this example the maximal nilpotent ideal is $N=\langle p \rangle$, with nilpotency index equal to $k$ in $\mathbb{Z}_{p^k}$, and $\F\cong \Z_p$. Thus, 
\begin{equation}\label{few1}
(\mathbb{Z}_{p^k}G)^{\ast}=\left\{g^{r}f^{r-1}: r=p^{k-1}, \bar{f},\bar{g}\in (\mathbb{Z}_{p} G)^{\ast} \hspace{0.1cm}\text{and}\hspace{0.1cm} \bar{f}\bar{g}=\bar{1}\right\}.
\end{equation}
If the group $G$ is cyclic of order $n$, it is known that $\mathbb{Z}_{p^k}G\cong \mathbb{Z}_{p^k}[x]/\langle x^n-1\rangle$, then (\ref{few1}) can be expressed in the following equivalent form,
\[
(\Z_{p^k}[x]/\langle x^n-1 \rangle)^{\ast}=\left\{g^{r}f^{r-1}: r=p^{k-1},\bar{f},\bar{g}\in (\Z_p[x]/\langle x^n-1 \rangle)^{\ast}\hspace{0.1cm}\text{and}\hspace{0.1cm} \bar{f}\bar{g}=\bar{1}\right\}.
\]
 
Another interesting example of chain rings is represented by Galois rings, $R=\mathbb{Z}_{p^k}[x]/\langle q(x) \rangle$ where $p$ is a prime number, $k>1$ and $q(x)$ is a monic polynomial of degree $r$ whose image in $\Z_{p}[x]$ is irreducible. In this example the maximal nilpotent ideal is $N=pR$, with nilpotency index equal to $k$ in $R$, and $\F=R/N\cong \F_{p^{r}}$. Thus, 
\begin{equation}\label{few2}
((\mathbb{Z}_{p^k}[x]/\langle q(x) \rangle)G)^{\ast}=\left\{g^{r}f^{r-1}: r=p^{k-1},\bar{f},\bar{g}\in (\F_{p^{r}}G)^{\ast}\right\}.
\end{equation}

\begin{example}
Consider the group ring $\Z_{25} C_5$, where $C_5$ is the cyclic group of order five generated by $a$. Observe that $N = \langle 5 \rangle$ is a nilpotent ideal of index $2$ in the ring $\Z_{25}$, so the collection $\{N C_5, N^2C_5\}$ of ideals of the ring $\Z_{25}C_5$ satisfies the CNC-codition with constant characteristic $s = 5$. Since the element $ f = 2 - a=2+24a$ is an invertible element in the group ring $\Z_{25} C_5$ and $g = 1 + 3a +4a^2 +2a^3 + a^4 \in (\Z_{25}/N) C_5$ is such that $\bar{f}\bar{g} = 1$ in $(\Z_{25}/N)C_5$, from Corollary \ref{GeralNilGR} it follows that $f + NC_5 \subset (\Z_{25} C_5)^{\ast}$ and $f^{-1} = g(fg)^4$. Thus, $fg=1+5a+5a^2$, $(fg)^4=1+20a+20a^2$ and, finally, 
\[f^{-1} = g(fg)^4 = 11+18a+9a^2+17a^3+21a^4\]
is the inverse of $f$ in the group ring  $\Z_{25} C_5$.

\end{example}


\subsection{Commutative group rings $\mathbb{Z}_mG$ with $G$ a commutative group}

The previously discussed results for $\Z_{p^k}G$ can be extended to the group ring $\Z_m G$, where $m>1$ and $G$ is a finite commutative group.

\begin{theo}\label{Zn}
Let $m=p_1^{r_1}p_2^{r_2}\cdots p_j^{r_j}$ be the prime factorization of the integer $m \geq 2$. 
Set $m_i=m/p_i^{r_i}$ and let $s_i$ be a natural number such that $s_im_i=1\mod(p_i^{r_i})$ for $i=1,2,3,\dots,j$. If $\bar f_i$ is an invertible element of the group ring $\mathbb{Z}_{p_i}G$ for $i=1,2,3,\dots,j$, then 
\begin{equation}\label{idew}
f=s_1m_1f_1+s_2m_2f_2+\cdots+s_jm_jf_j,
\end{equation}
is an invertible element in the group ring $\mathbb{Z}_{m}G$, with 
\begin{equation}\label{ide}
f^{-1}=s_1m_1g_1^{\alpha_1}f_1^{\alpha_1-1}+s_2m_2g_2^{\alpha_2}f_2^{\alpha_2-1}+\cdots+s_jm_jg_j^{\alpha_j}f_j^{\alpha_j-1},
\end{equation}
where $\alpha_i=p_i^{r_i-1}$ and $\bar g_i\in \mathbb{Z}_{p_i}G$ is such that $\bar f_i\bar g_i=\bar 1$. Moreover,
\begin{equation}\label{euip}
|(\mathbb{Z}_{m}G)^{\ast}|=(m/(p_1p_2\cdots p_j))^{|G|}|(\mathbb{Z}_{p_1}G)^{\ast}||(\mathbb{Z}_{p_2}G)^{\ast}|\cdots|(\mathbb{Z}_{p_j}G)^{\ast}|.
\end{equation}
\end{theo}
\proof From the Chinese Remainder Theorem (CRT), $\Z_m\cong \Z_{p_1^{r_1}}\times \Z_{p_2^{r_2}}\times \cdots\times\Z_{p_j^{r_j}}$, and 
\begin{equation*}\label{is}
\mathbb{Z}_{m}G\cong_{\phi}\mathbb{Z}_{p_1^{r_1}}G \times \mathbb{Z}_{p_2^{r_2}}G \times\cdots\times \mathbb{Z}_{p_j^{r_j}}G.
\end{equation*}
Therefore,
\begin{equation}\label{is}
(\mathbb{Z}_{m}G)^{\ast}\cong_{\phi}(\mathbb{Z}_{p_1^{r_1}}G)^{\ast} \times (\mathbb{Z}_{p_2^{r_2}}G)^{\ast} \times\cdots\times(\mathbb{Z}_{p_j^{r_j}}G)^{\ast}.
\end{equation}
If $\bar{f_i}$ is an invertible element in $\mathbb{Z}_{p_i}G$, from relation (\ref{few1}) $f_i$ is an invertible element in the group ring $\mathbb{Z}_{p_i^{r_i}}G$  with $f_i^{-1}=g_i^{\alpha_i}f_i^{\alpha_i-1}$, where $\bar f_i\bar g_i=\bar 1$ for each $i=1,2,3,\cdots,j$. Thus,  $h=(f_1,f_2,\cdots,f_j)$  is an invertible element in the product group given in (\ref{is}), with 
\[
h^{-1}=(g_1^{\alpha_1}f_1^{\alpha_1-1},g_2^{\alpha_2}f_2^{\alpha_2-1},\cdots,g_j^{\alpha_j}f_j^{\alpha_j-1}).
\]
Consequently, $f=\phi^{-1}(h)$ is an invertible element in the group ring $\mathbb{Z}_{m}G$, with $f^{-1}=\phi^{-1}\left(h^{-1}\right)$. Finally, from the CRT, $f$ and $f^{-1}$ can be expressed in the form (\ref{idew}) and (\ref{ide}) respectively. The equality in (\ref{euip}) follows from (\ref{fe}) and (\ref{is}).  
\cqd

The following result provides an alternative way to compute the invertible elements in $\Z_mG$.

\begin{theo}\label{Znv2}
Let $m=p_1^{r_1}p_2^{r_2}\cdots p_j^{r_j}$ be the prime factorization of the integer $m \geq 2$. 
Set $k=\max\{r_1,r_2,\ldots,r_j\}$, $c_i=(p_1p_2\cdots p_j)/p_i$, and let $t_i$ be a natural number such that $t_ic_i=1\mod(p_i)$ for $i=1,2,3,\cdots,j$. If $\bar f_i$ is an invertible element of the group ring $\mathbb{Z}_{p_i}G$ for $i=1,2,3,\dots,j$, then 
\begin{equation}\label{idew1}
f=t_1c_1 f_1+t_2c_2 f_2+\cdots+t_jc_j f_j,
\end{equation}
is an invertible element in the group ring $\mathbb{Z}_{m}G$, with
\begin{equation}\label{ide1}
f^{-1}=(t_1c_1g_1+t_2c_2g_2+\cdots+t_jc_jg_j)^{(p_1p_2\cdots p_j)^{k-1}}(t_1c_1f_1+t_2c_2f_2+\cdots+t_jc_jf_j)^{(p_1p_2\cdots p_j)^{k-1}-1} ,
\end{equation}
where $\bar g_i\in \mathbb{Z}_{p_i}G$ is such that $\bar f_i\bar g_i=\bar 1$.
Moreover,
\begin{equation}\label{equipv2}
|(\mathbb{Z}_{m}G)^{\ast}|=(m/(p_1p_2\cdots p_j))^{|G|}|(\mathbb{Z}_{p_1p_2\cdots p_j}G)^{\ast}|.
\end{equation}
\end{theo}
\proof If $\bar{f_i}$ is an invertible element in the group ring $\mathbb{Z}_{p_i}G$ for $i=1,2,3,\cdots,j$, from Theorem \ref{Zn} it follows that 
\begin{equation}\label{idfg}
f=t_1c_1f_1+t_2c_2f_2+\cdots+t_jc_jf_j
\end{equation}
is an invertible element in the group ring $\mathbb{Z}_{p_1p_2\cdots p_j}G$, with
\begin{equation}\label{idfg1}
w:=f^{-1}=t_1c_1g_1+t_2c_2g_2+\cdots+t_jc_jg_j,
\end{equation}
where $\bar g_i\in \mathbb{Z}_{p_i}G$ is such that $\bar f_i\bar g_i=\bar 1$. By observing that $a=p_1p_2\cdots p_j$ has nilpotency index $k$ in the ring $\mathbb{Z}_m$ and,
\begin{equation}\label{ort}
\frac{\mathbb{Z}_m}{a\mathbb{Z}_m}\cong \mathbb{Z}_a
\end{equation}
has characteristic $a=p_1p_2\cdots p_j$, from Corollary \ref{cyu} it follows that $f$ is an invertible element in $\mathbb{Z}_{m}G$, with
\[ \label{idev2}
f^{-1}=w^{(p_1p_2\cdots p_j)^{k-1}}f^{(p_1p_2\cdots p_j)^{k-1}-1}.
\] Relation (\ref{equipv2}) follows from Corollary \ref{cyu} and (\ref{ort}).
\cqd


\section*{Acknowledgements}
This research was carried on while the first author was a post-doctoral scholar at the Departamento de Matem\'aticas, Universidad Aut\'onoma Metropolitana - Iztapalapa, Cd. de M\'exico, M\'exico. The second author was a visiting researcher at the Instituto de Investigaciones en Matem\'aticas Aplicadas y Sistemas (IIMAS), Universidad Nacional Aut\'onoma de M\'exico, Cd. de M\'exico. They would like to thank these institutions for the warm reception during their stay.



\begin{thebibliography}{19}

\bibitem{ch-l} Chebolu, S. K. and Lockridge, K., ``Fuchs' problem for $p$-groups", arXiv: 1901.10081v1 [Math-RA], 29 Jan 2019.

\bibitem{con} Conrand, K., ``Dirichlet's Unit Theorem", Retrieved from: \\ https://kconrad.math.uconn.edu/blurbs/gradnumthy/unittheorem.pdf

\bibitem{delc-d} Del Corso, I. and Dvornicich, R., ``On Fuchs' problem about the group of units of a ring",  {\it Bull. London Math Soc.} 50, (2018) 274-292

\bibitem{dolz} Dolzan, D., ``Group of units in a finite ring", {\it J. of Pure and Applied Algebra} 170, (2002), 175-183.

\bibitem{h-l-m} Hou, X-D, Leung, K.H. and Ma, S.L. ``On the group of units of finite commutative chain rings", {\it Finite Fields and Applic.} No.9, (2003) 20-38.

\bibitem{jacobson} Jacobson, N., \textit{Basic  Algebra II}. W. H. Freeman and Company, New York, (1989).

\bibitem{k-k-s} Kanwar, P., Khatkar, M., and Sharma, R.K., ``Idempotent and units of matrix rings over polynomial rings", {\it Int. Elec. Jour. of Algebra}, 22(2017) 147-169. 

\bibitem{k} Kline, D., ``The structure of unit groups" Retrieved from: \\ http://math.uchicago.edu/~may/REU2014/REUPapers/Kline.pdf.

\bibitem{mcdonald} McDonald,B. R., {\it Finite Rings with Identity.} Marcel Dekker, Inc., New York, 1974.

\bibitem{liftidemp} de Melo Hern\'andez, F.D., Hern\'andez Melo, C.A. and Tapia-Recillas, H., ``On idempotents of a class of commutative rings", {\it Communications in Algebra}, (2020)  DOI: 10.1080/00927872.2020.1754424.

\bibitem{m-m} Menal, P. and Moncasi, J., ``Lifting units in self-injective rings and index theory for Rickart $C^{\ast}$-algebras". {\it Pacific Journal of Math.}, 126, No.2, (1987) 295-329.

\bibitem{perera} Perera, F., ``Lifting units modulo exchange ideals and $C^{\ast}$-algebras with real rank zero", arXiv: 9906181v2 [math.RA] 25 Aug 1999.

\bibitem{janez} Ster, J., ``Lifting units in clean rings", {\it Journal of Algebra}, 381, (2013) 200-208,


\end{thebibliography}
\end{document}